\documentclass[11pt,reqno]{amsart}
\usepackage{mathtools}
\usepackage{xcolor}
\usepackage{amsmath}
\usepackage{amsthm}
\usepackage{amsopn}
\usepackage{amssymb}
\usepackage{fullpage}
\usepackage{enumerate}
\numberwithin{equation}{section}
\usepackage{float}
\usepackage{graphicx}
\usepackage{titlesec}
\usepackage{fancyhdr}
\usepackage{hyperref}

\usepackage{soul}
\usepackage{cancel}

\newtheorem{thm}{Theorem}

\newtheorem{theorem}{Theorem}[section]

\newtheorem{lemma}[theorem]{Lemma}

\theoremstyle{definition}
\newtheorem{definition}[theorem]{Definition}

\DeclareMathOperator{\Capac}{cap}

\DeclareMathOperator{\T}{\mathbb{T}}

\DeclareMathOperator{\D}{\mathbb{D}}

\DeclareMathOperator{\Hol}{Hol}

\titleformat{\section}
  {\normalfont\scshape}{\thesection}{1em}{\centering}
\titleformat{\subsection}[runin]
  {\bfseries}{\thesubsection}{1em}{}
\DeclareMathOperator{\N}{\mathbb{N}}

\DeclareMathOperator{\C}{\mathbb{C}}
\newcommand{\abs}[1]{\left| #1 \right|}

\title{Hardy number of Koenigs domains: sharp estimate}
\date{\today}

\author[M.D. Contreras]{Manuel D. Contreras}
\address{Manuel D. Contreras. Departamento de Matem\'atica Aplicada II and IMUS, Escuela T\'ecnica Superior de Ingenier\'ia, Universidad de Sevilla, Camino de los Descubrimientos, s/n 41092, Sevilla, Spain}
\email{contreras@us.es}  

\author[F. J. Cruz-Zamorano]{Francisco J. Cruz-Zamorano}
\address{Francisco J. Cruz-Zamorano. Departamento de Matem\'atica Aplicada II and IMUS, Escuela T\'ecnica Superior de Ingenier\'ia, Universidad de Sevilla, Camino de los Descubrimientos, s/n 41092, Sevilla, Spain}
\email{fcruz4@us.es}  

\author[M. Kourou]{Maria Kourou}
\address{Maria Kourou. Julius-Maximilians-Universit\"at W\"urzburg, Institut f\"ur Mathematik, Emil Fischer Stra{\ss}e 40, 97074, W\"urzburg, Germany}
\email{maria.kourou@uni-wuerzburg.de}   

\author[L. Rodr\'iguez-Piazza]{Luis Rodr\'iguez-Piazza}
\address{Luis Rodr\'iguez-Piazza. Departmento de An\'alisis Matem\'atico and IMUS, Facultad de Matem\'aticas, Universidad de Sevilla, Calle Tarfia, s/n 41012 Sevilla, Spain}
\email{piazza@us.es}  

\thanks{M. D. Contreras, F. J. Cruz-Zamorano and L. Rodr\'iguez-Piazza are partially supported by Ministerio de Innovaci\'on y Ciencia, Spain, project PID2022-136320NB-I00. }
\thanks{F. J. Cruz-Zamorano is also partially supported by Ministerio de Universidades, Spain, through the action Ayuda del Programa de Formaci\'on de Profesorado Universitario, reference FPU21/00258.}
\thanks{M. Kourou is partially supported by the Alexander von Humboldt Foundation.}

\newtheorem{notation}[theorem]{Notation}
\DeclareMathOperator{\h}{h}

\begin{document}
\begin{abstract}
Let $\Omega$ be a regular Koenigs domain in the complex plane $\C$. We prove that the Hardy number of $\Omega$ is greater or equal to $1/2$. That is, every holomorphic function in the unit disc $f:\D\to \Omega$ belongs to the Hardy space $H^{p}(\D)$ for all $p<1/2$. 
\end{abstract}
\maketitle

\section{Introduction}

Given a holomorphic function in the unit disc $f:\D\to \C$
one can define its \textit{Hardy number}  by
$$\h(f) := \sup\left(\{ p>0 : f \in H^p(\D) \} \cup \{0\}\right) \in [0,+\infty],$$
which somehow measures the ``growth'' of $f$. As usual, given $p\in (0,+\infty]$, $H^p(\D)$ denotes the Hardy space. Motivated by this concept, Hansen \cite{Hansen1970} introduced a similar idea in order to examine the range of holomorphic functions taking values in a domain $\Omega$. The so-called \textit{Hardy number} of a domain $\Omega$ is defined as
$$\h(\Omega) := \inf \{\h(f):  f\in \Hol(\D,\Omega)  \}.$$

It is worth pointing out that in \cite{CCKR}, it was shown that for every  $p\in  [0,+\infty]$ there exists a domain $\Omega $ in the complex plane such that $\h(\Omega)=p$. 

A classical problem in the theory of Hardy spaces has been to determine the range of the Hardy number of a hyperbolic planar domain $\Omega$ (i.e., those whose boundary contains at least two points) in terms of its geometry and boundary behavior. As a matter of fact, several works have contributed to estimates on the Hardy number for certain classes of hyperbolic domains. Hansen \cite{Hansen1970, Hansen1971} provided a characterization of the Hardy number for starlike and spirallike with respect to the origin domains. A few years later, Ess\'en \cite{Essen_1981} obtained estimates of the Hardy number of a general hyperbolic domain in terms of harmonic measures and the logarithmic capacity. Thereafter, Kim and Sugawa \cite{Kim_Sugawa2011} examined the Hardy number for unbounded $K$-quasidiscs. 
Quite recently, the range of the Hardy number for comb domains was studied by Karafyllia \cite{Kar_2022}.

Inspired by the impact of Schr\"oder's  functional equation on the discrete dynamics of the unit disc,
 Poggi-Corradini \cite{PC_1997_2} examined the Hardy number of hyperbolic domains $\Omega$ satisfying Schr\"oder's inclusion property; i.e. $\lambda \Omega \subseteq \Omega$, for some $\lambda \in \D$. 
 In a complementary context in holomorphic dynamics in the unit disc, one finds the Abel's functional equation. Driven by this equation, in the  very recent paper \cite{CCKR}, it was analyzed the Hardy number of hyperbolic domains $\Omega$ satisfying \textit{Abel's inclusion property}; namely $\Omega+1\subseteq \Omega$. From now on, these domains will be called \textit{Koenigs domains}. Among other results, in that paper it was shown that there exists $C > 0$ such that, for all regular Koenigs domains $\Omega$, one has $\h(\Omega) \geq C$ (see \cite[Theorem 1.1]{CCKR}). 
 Let us recall that a domain $\Omega$ in $\C$ is called {\sl regular} if the logarithmic capacity of $\C\setminus \Omega$ is positive.  
 Using elementary examples of Hardy numbers, it is clear that $C\leq 1/2$. In addition, it is well-known that the Hardy number of a non-regular domain is zero. The question of determining the sharp value of the constant $C$ arises naturally. Our main result answers this problem providing that, sure enough, $C=1/2$:
 
\begin{theorem}
\label{thm:main-hardy}
Let $\Omega \subset \C$ be a Koenigs domain. Then, $\h(\Omega) \geq 1/2$ if and only if $\Omega$ is a regular domain. In particular, for every Koenigs domain $\Omega$ it holds that $\h(\Omega)\in \{0\}\cup [1/2,+\infty]$. 
\end{theorem}

The proof of this result strongly depends on potential theory, namely a sharp estimate of the logarithmic capacity of compact sets obtained as  union of integer translations of a compact non-regular  subset of $\C$. Recall that a compact set in the complex plane is {\sl polar} if its logarithmic capacity is zero.

\begin{theorem}
\label{thm:main-potential}
Let $E$ be a compact non-polar subset of $\overline{D(0,1/4)}$. Set $$K_n = \bigcup_{j=1}^n(E+j).$$
Then, 
$$\lim_{n \to +\infty}\dfrac{\Capac(K_n)}{\Capac([0,n])} = 1.$$
Moreover, $\log(\Capac(K_n)) = \log(n/4) + \mathcal{O}(1/\sqrt{n})$.
\end{theorem}

The structure of the paper is as follows: after a preliminary introduction on logarithmic capacity, Hardy number and harmonic measure in Section \ref{sec:preliminaries}, we prove Theorem \ref{thm:main-potential} in Section \ref{sec:proof1}.
Finally, in Section \ref{sec:proof2}, we prove Theorem \ref{thm:main-hardy}.

\begin{notation}
We write $a_n = \mathcal{O}(b_n)$ if there exists $C > 0$ and $N \in \N$ such that $\abs{a_n} \leq Cb_n$ for all $n \geq N$.
\end{notation}

\section{Preliminaries}\label{sec:preliminaries}
\subsection{Logarithmic potentials.}
In Potential Theory, one of the main tools are the potentials of measures, which can be defined as follows:
\begin{definition}
\label{def:potential-energy}
Let $\mu$ be a finite positive measure on $\C$ with compact support. Its (logarithmic) potential is the function $p_{\mu} \colon \C \to [-\infty,+\infty)$ given by
$$p_{\mu}(z) = \int_{\C}\log\abs{z-w}d\mu(w), \quad z \in \C.$$
Its (logarithmic) energy $I(\mu) \in [-\infty,+\infty)$ is given by
$$I(\mu) = \int_{\C}p_{\mu}(z)d\mu(z) = \int_{\C}\int_{\C}\log\abs{z-w}d\mu(w)d\mu(z).$$
From these concepts one can define the (logarithmic) \textit{capacity} of a set $X \subseteq \C$, that is
\begin{equation}\label{eq:capacitydef}
    \Capac(X) = \sup_{\mu}\exp(I(\mu))
\end{equation}
where the supremum is taken among every probability measure $\mu$ whose support is a compact subset of $X$.
\end{definition}

In these definitions, we understand that $I(\mu) = -\infty$ if the former integral is not convergent. Indeed, if $\Capac(X) = 0$, then $I(\mu) = -\infty$ for every probability measure $\nu$ with compact support in $X$. If this happens, $X$ is said to be a \textit{polar set}.
A property is said to be satisfied \textit{nearly everywhere} in $X \subset \mathbb{C}$, if it is satisfied for all points in $X$, except maybe for a Borel polar subset. In general, polar sets are negligible from the potential-theoretic point of view. Indeed, the following lemma holds:
\begin{lemma}
\cite[Corollary 3.2.5]{Ransford}
\label{lem:countable}
A countable union of Borel polar sets is polar.
\end{lemma}

Let $\Omega \subseteq \C$ be a domain, that is, an open and connected set. If $\Capac(\C \setminus \Omega) > 0$, then $\Omega$ is called \textit{regular}. Furthermore, when it comes to compacta, the theory of the potentials is richer: if $X\subset \C$ is a compact set, then there exists a probability measure $\nu$ with support in $X$ attaining the supremum in \eqref{eq:capacitydef}. In fact, this measure is unique if $X$ is non-polar and in such a case $\nu$ is said to be the \textit{equilibrium measure} of $X$; see \cite[Theorem 3.7.6]{Ransford}. The potential associated to the equilibrium measure satisfies the following property:
\begin{thm}[Frostman's Theorem]
\cite[Theorem 3.3.4]{Ransford}
\label{thm:Frostman}
Let $X \subset \C$ be a non-polar compact set, and let $\nu$ be its equilibrium measure. Then, $p_{\nu}(z) \geq I(\nu)$ for all $z \in \C$. Moreover, $p_{\nu} \equiv I(\nu)$ nearly everywhere on $X$.
\end{thm}

\subsection{Hardy numbers.}
\label{subsec:HardyNumbers}
Given $0 < p < +\infty$, the Hardy space $H^p(\D)$ is defined as the set of all holomorphic maps $f \colon \D \to \C$ such that
$$\sup_{0<r<1}\int_{\T}\abs{f(r\xi)}^pdm(\xi) < +\infty,$$
where $m$ is the normalized length measure in the boundary of the unit disc $\T$. In the case $p = +\infty$, the Hardy space $H^{\infty}(\D)$ stands for the set of all holomorphic maps $f \colon \D \to \C$ that are bounded, that is,
$$\sup_{z \in \D}\abs{f(z)} < +\infty.$$
We refer to \cite{Duren_Hp} for a complete introduction to this topic. One of the properties of these spaces is that they form a decreasing family, that is, $H^p(\D) \supseteq H^q(\D)$ if $0 < p \leq q \leq +\infty$. These relations suggest the following idea: given a holomorphic map $f \colon \D \to \C$, its Hardy number $\h(f)$ is defined as
$$\h(f) = \sup(\{0\} \cup \{p > 0 : f \in H^p(\D)\}) \in [0,+\infty].$$
Note that $f \in H^p(\D)$ for every $0 < p < \h(f)$, and $f \not\in H^p(\D)$ if $p > \h(f)$.

In a similar manner, this idea can be translated to domains. Given a domain $\Omega \subset \C$, its Hardy number is defined as
$$\h(\Omega) = \inf\{\h(f) : f \in \text{Hol}(\D,\Omega)\}.$$
Here, $\text{Hol}(\D,\Omega)$ denotes the set of all holomorphic maps $f \colon \D \to \C$ such that $f(\D) \subseteq \Omega$. Let us state some well-known properties of the Hardy number of a domain:
\begin{lemma}
\cite[Lemmas 2.1 and 2.3]{Kim_Sugawa2011}
\label{lem:properties}
Let $\Omega, \Omega^{\prime} \subseteq \C$ be two domains. Then 
\begin{enumerate}[\normalfont(a)]
\item $\h(\Omega^{\prime}) \leq \h(\Omega)$, if $\Omega\subseteq \Omega^{\prime}$.
\item $\h(\varphi(\Omega))=\h(\Omega)$ for a complex affine map $\varphi(z)=az+b$, $a\neq 0$. 
\item $\h(\Omega) \geq \frac{1}{2}$ if $\Omega \neq \C$ is simply connected.
\item $\h(\Omega) = \h(p)$, where $p \colon \D \to \Omega$ is a universal covering map of $\Omega$.
\end{enumerate}
\end{lemma}

A characterization of the Hardy number for certain domains can be obtained through a classical result from Nevanlinna.
\begin{thm}
\label{thm:Nevanlinna}
\cite[Theorems 5.1.1 and 5.4.2, p. 209 and 211]{Nevanlinna1970}
A universal covering map $p \colon \D \to \Omega$ has non-tangential limits almost everywhere on $\T$ if and only if $\Omega$ is a regular domain.
\end{thm}
Recall that any function whose Hardy number is positive (that is, it belongs to some Hardy space) has non-tangential limits almost everywhere on $\T$. Then, from the above result and the properties of the Hardy numbers, it follows that the Hardy number of every non-regular domain is zero.

\subsection{Harmonic measure.}
For a regular domain $\Omega \subseteq \C$, the harmonic measure $\omega(z,B,\Omega)$ of a Borel set $B \subseteq \partial\Omega$ with respect to a point $z \in \Omega$ is the value at $z$ of the solution of the generalized Dirichlet problem in $\Omega$ with boundary values $1$ on $B$ and $0$ on $\partial \Omega \setminus B$. The name comes after the following properties:
\begin{enumerate}[\normalfont(i)]
\item For a fixed Borel set $B \subseteq\partial \Omega$, the map $z \mapsto \omega (z, B, \Omega)$ is harmonic and bounded.
\item For a fixed point $z \in \Omega$, the map $B \mapsto \omega(z, B, \Omega)$ is a probability measure on $\partial \Omega$.
\end{enumerate}
For an introduction to harmonic measure, we refer to \cite[Section 4.3]{Ransford}.

A relation between the Hardy number of a domain and the limiting behaviour of certain harmonic measures was proposed by Ess\'en in \cite[Lemma 1]{Essen_1981}. Following his ideas, Kim and Sugawa proved the following useful formula for the Hardy number:
\begin{thm}
\label{thm:KimSugawa}
\cite[Lemma 3.2]{Kim_Sugawa2011}
Let $\Omega \subseteq\C$ be a domain with $0 \in \Omega$. Then,
$$\h(\Omega) = \liminf_{R \to + \infty}\left(-\dfrac{\log \omega(0,F_R,\Omega_R)}{\log R}\right),$$
where $\Omega_R$ is the connected component of $\Omega \cap D(0,R)$ containing the origin and $F_R = \partial\Omega_R \cap \{w \in \C : |w| = R\}$.
\end{thm}

\subsection{Complex dynamics on the unit disc.}
A function $\phi \colon \D \to \D$ is called \textit{non-elliptic} if it has no fixed points, that is, $\phi(z) \neq z$ for every $z \in \D$. For a complete picture of the dynamical aspects of holomorphic maps on the unit disc, we refer to \cite{AbateBook}. For example, the following result relates non-elliptic maps and translations:
\begin{theorem}
\cite[Theorem 4.6.8]{AbateBook}
For every non-elliptic map $\phi$ there exists a map $\sigma \colon \D \to \C$ that semiconjugates $\phi$ with a unit translation, that is, $\sigma \circ \phi = \phi + 1$.
\end{theorem}
The latter theorem merges some important contributions of Valiron \cite{Valiron} and Pommerenke \cite{Pom_1979}, who also collaborated with Baker \cite{BakerPommerenke}, over the last century. Independently analyzing different cases, they all showed iterative constructions for the map $\sigma$, which is known as a Koenigs function for $\phi$. 

For a Koenigs function $\sigma$, its image domain $\Omega = \sigma(\D)$ satisfies an important inclusion property. If $w = \sigma(z) \in \Omega$ for some $z \in \D$, then $w+1 = \sigma(z)+1 = \sigma(\phi(z)) \in \Omega$. That is, $\Omega + 1 \subseteq \Omega$. Thus, we give the following definition:
\begin{definition}
A domain $\Omega \subseteq \C$ is said to be a Koenigs domain if $\Omega + 1 \subseteq \Omega$.
\end{definition}
Notice that $\h(\sigma) \geq \h(\sigma(\D))$. Therefore, it is possible to analyze the membership of a Koenigs function to some Hardy space by analyzing the Hardy number of its image Koenigs domain.

\section{Proof of Theorem \ref{thm:main-potential}}\label{sec:proof1}

Before moving on to the proof of Theorem \ref{thm:main-potential}, we prove two auxiliary lemmas related to the construction of the equilibrium measure of the compact sets $K_n$, $n \in \mathbb{N}$.

In the sequel we will use the equilibrium measure of a compact interval, which can easily be derived from \cite[Eq. (1.7), p. 25]{SaffTotik}. For any $n \in \N$, the equilibrium measure $\mu$ of the interval $[0,n]$ is absolutely continuous (with respect to Lebesgue's measure $m$) and it is given by
\begin{equation}
\label{eq:mu}
\dfrac{d\mu}{dm}(t) = \dfrac{\chi_{[0,n]}(t)}{\pi \sqrt{t(n-t)}}.
\end{equation}
In particular, it follows
\begin{equation}
\label{eq:properties_mu}
p_{\mu}(t) = I(\mu) = \log(\Capac([0,n])) = \log(n/4), \quad \text{for all } t \in [0,n]\, ; 
\end{equation}
the above result is a combination of Theorem \ref{thm:Frostman} and \cite[Theorems 4.2.2 and 4.2.4]{Ransford}.

Using this equilibrium measure, we define the following coefficients:
$$\alpha_{j}:= \int_{j-1}^jd\mu(t) = \frac{1}{\pi} \int_{j-1}^{j} \frac{dt}{\sqrt{t(n-t)}}, \quad j \in \N, \quad 1 \leq j \leq n.$$
Notice that $\alpha_j$ also depends on $n$. However, in seek of clearance, this it not explicitly written in the notation. Some properties of these numbers can be derived directly from the definition. For example, it follows that $\alpha_j>0$ for all $j=1,..., n$ and that $\sum_{j=1}^n \alpha_{j}=1$. It is also possible to notice that these numbers are endowed with some symmetry, namely $\alpha_j = \alpha_{n-j+1}$, and that $\alpha_j$ is non-increasing for $j=1, ..., \lfloor (n+1)/2 \rfloor$, where $\lfloor x \rfloor$ denotes the integer part of the real number $x$.

We prove the following estimations:

\begin{lemma}\label{lem:properties_aj}
\begin{enumerate}[\normalfont(a)]
\item There exists $C_1 > 0$ such that $\alpha_j \leq C_1/\sqrt{n}$, for all $n \in \N$ and all $j = 1,\ldots,n$.
\item There exists $C_2 > 0$ such that $$\sum_{\substack{j=1 \\ j\neq k}}^n \frac{\alpha_j}{|j-k|} \leq \frac{C_2}{\sqrt{n}},$$ for all $n \in \N$ and all $k = 1,\ldots,n$.

\end{enumerate}
\begin{proof}
(a) Fix a natural number $n\geq 2$. Then 
$$\alpha_1 = \dfrac{1}{\pi} \int_0^1\dfrac{dt}{\sqrt{t(n-t)}} \leq \dfrac{1}{\pi\sqrt{n-1}}\int_0^1\dfrac{dt}{\sqrt{t}} = \mathcal{O}\left(\dfrac{1}{\sqrt{n}}\right).$$
Thanks to the monotonicity and symmetric properties of the coefficients $\alpha_j$, we have that $\alpha_j\leq \alpha_1$, for $j=1, ...,n$, so that (a) holds.

(b) Fix $n \in \N$ and define
$$S(k) = \sum_{\substack{j=1 \\ j\neq k}}^n \frac{\alpha_j}{|j-k|}, \quad k = 1,\ldots,n.$$
By symmetry, notice that $S(k) = S(n-k+1)$. Thus, it is enough to work with $k \in \N$ such that $k \leq (n+1)/2$. If this is the case, notice that for any $j \in \N$ with $1 \leq j \leq (n+1)/2$ it is possible to check that $\alpha_j = \alpha_{n-j+1}$ but $\abs{j-k} \leq \abs{n-j+1-k}$. Therefore,
$$S(k) \leq 2 \sum_{\substack{j=1 \\ j\neq k}}^{\lfloor (n+1)/2 \rfloor} \frac{\alpha_j}{|j-k|}.$$

Given two finite sequences $\{ a_1, a_2, ..., a_m \}$ and $\{ b_1, b_2, ..., b_m \}$ of non-negative real numbers, the Hardy-Littlewood inequality asserts that
$$
\sum _{j=1}^m a_j b_j\leq \sum _{j=1}^m a_j^* b_j^*,
$$
where $\{a_j^*\}$ denote the sequence of elements $a_j$ arranged in decreasing order and similarly for  $\{b_j^*\}$ (see \cite[\S 10.2]{inequalities} and also \cite[p. 43]{Bennett-Colin}).
Notice that the map $\{1,\ldots,\lfloor (n+1)/2 \rfloor\} \ni j \mapsto \alpha_j$ is decreasing. Moreover, taking $b_j = 1/|j-k|$, for $j=1, ...,  \lfloor (n+1)/2 \rfloor$, $j\neq k$, and $b_{k}=0$, we have that $b_j^*\leq 2/j$. Thus
$$\sum_{\substack{j=1 \\ j\neq k}}^{\lfloor (n+1)/2 \rfloor} \frac{\alpha_j}{|j-k|} =
 \sum_{j=1}^{\lfloor (n+1)/2 \rfloor} \alpha_j b_{j} \leq 2\sum_{j=1}^{\lfloor (n+1)/2 \rfloor}\dfrac{\alpha_j}{j}.$$

But now, notice that, for $n\geq 2$,
\begin{equation*}
\begin{split} \sum_{j=1}^{\lfloor (n+1)/2 \rfloor}\dfrac{\alpha_j}{j}&  \leq 
\dfrac{2}{\pi}\int_0^{(n+1)/2}\dfrac{dt}{(t+1)\sqrt{t(n-t)}} \leq \dfrac{2}{\pi\sqrt{(n-1)/2}}\int_0^{+\infty}\dfrac{dt}{(t+1)\sqrt{t}} = \mathcal{O}\left(\dfrac{1}{\sqrt{n}}\right), 
\end{split}
\end{equation*}
where, in the first inequality, we have used that $1/j \leq 2/(t+1)$ whenever $j-1\leq t\leq j$. Thus, the result follows.
\end{proof}
\end{lemma}

To proceed with the proof of Theorem \ref{thm:main-potential}, let $E \subseteq \overline{D(0,1/4)}$ be a compact non-polar set with $0\in E$ and equilibrium measure $\nu$. Fix $n \in \mathbb{N}$ and set
$$K_n=\bigcup_{j=1}^n (E+j).$$
Let us define the positive measure $\sigma$ given by
$$\sigma(A) = \sum_{j=1}^{n} \alpha_j \nu(A-j),$$
where $A \subset \C$ is a Borel set. Notice that $\sigma$ is a probability measure whose support is a compact set lying on $K_n$.
Once more, observe that $\sigma$ depends on $n$.

Recalling the Definition \ref{def:potential-energy}, let us prove the following lemma:
\begin{lemma}\label{lemma:sigma}
Under the above notation, the following statements hold:
\begin{enumerate}[\normalfont(a)]
\item $p_{\sigma}(x) \geq \log(n/4) + \mathcal{O}(1/\sqrt{n})$ for every $x \in K_{n}$.
\item $p_{\sigma}(x) = \log(n/4) + \mathcal{O}(1/\sqrt{n})$ for nearly every $x \in K_n$.
\item $\abs{p_{\sigma}(x)-p_{\sigma}(y)} = \mathcal{O}(1/\sqrt{n})$ for nearly every $x,y \in K_n$.
\item $p_{\sigma}(x) = \log(n/4) + \mathcal{O}(1/\sqrt{n})$ for every $x \in E$.
\item $\abs{p_{\sigma}(x)-p_{\sigma}(y)} = \mathcal{O}(1/\sqrt{n})$ for every $x \in E$ and nearly every $y \in K_n$.
\end{enumerate} 
In fact, the underlying constants do not depend on $x$ and $y$.
\begin{proof} In this proof we will always assume that $n\geq 2$. Note that, by the definition of the probability measure $\sigma$, its potential $p_{\sigma}$ can be written as
$$p_{\sigma}(z) = \sum_{j=1}^n  \alpha_j p_{\nu}(z-j) \, , \quad z \in \C.$$

Let $k \in \{1,...,n\}$. Notice that, by Theorem \ref{thm:Frostman},  $p_{\nu}(x-k) \geq  I(\nu)$ for all $x\in E+k$ and $p_{\nu}(x-k) =  I(\nu)$ nearly everywhere in $x\in E+k$. Therefore \begin{align*}
p_{\sigma}(x) & = \sum_{j=1}^n \alpha_j p_{\nu}(x-j) \geq \alpha_k I(\nu) + \sum_{\substack{j=1 \\ j\neq k}}^n \alpha_j \int_E \log|x-j-y| d\nu(y)
\end{align*}
for all $x\in E+k$ and the inequality turns to be an equality  nearly everywhere. We claim that, given $x\in E+k$, 
\begin{equation} \label{eq:claim}
\alpha_k I(\nu) + \sum_{\substack{j=1 \\ j\neq k}}^n \alpha_j \int_E \log|x-j-y| d\nu(y)= \log\left(\dfrac{n}{4}\right) + \mathcal{O}\left(\dfrac{1}{\sqrt{n}}\right).
\end{equation}
Thus, assuming this claim,  we clearly conclude (a) and (b).

Let us prove the claim. Take $x\in E+k$ and write $x^{\prime} =x-k \in E$. Then
\begin{align*}
\log|x-j-y| & = \log\abs{x^{\prime} -y + k-j} = \log\abs{k-j} + \log\abs{1+ \frac{x^{\prime} -y}{k-j}} \\
& = \log\abs{k-j} + \mathcal{O}\left(\frac{1}{\abs{k-j}} \right),
\end{align*}

for all $y \in E$, since $\abs{x^{\prime}-y} \leq 1/2$. 
Due to Lemma \ref{lem:properties_aj}.(b), we obtain 
\begin{equation}\label{eq:psigma}
\alpha_k I(\nu) + \sum_{\substack{j=1 \\ j\neq k}}^n \alpha_j \int_E \log|x-j-y| d\nu(y) = \alpha_k I(\nu) + \sum_{\substack{j=1 \\ j \neq k}}^n\alpha_j  \log|k-j| + \mathcal{O} \left(\frac{1}{\sqrt{n}} \right) \,  .
\end{equation}

Using the equilibrium measure $\mu$ for $[0,n]$ as described in \eqref{eq:mu}, let us define
\begin{align*}
S & := \alpha_k I(\nu) + \sum_{\substack{j=1 \\ j \neq k}}^n  \alpha_j\log|k-j|  - p_{\mu}(k-1/2) \\
& =  \alpha_k I(\nu) + \sum_{\substack{j=1 \\ j \neq k}}^n \alpha_j\log|k-j|  - \sum_{j=1}^n \int_{j-1}^j \log \left| k-1/2-t \right| d\mu(t) \\
& = \alpha_k I(\nu) - \int_{k-1}^k \log \left|k-1/2-t \right|d\mu(t) + \sum_{\substack{j=1 \\ j \neq k}}^n  \int_{j-1}^j  \left(\log\abs{k-j} - \log \left|k-1/2-t \right|\right)d\mu(t).
\end{align*}
We want to prove  that $\abs{S} = \mathcal{O}(1/\sqrt{n})$. This can be done by examining each of the three terms that define $S$. By Lemma \ref{lem:properties_aj}.(a), since $\abs{I(\nu)} < +\infty$, one concludes that $\abs{\alpha_k I(\nu)} = \mathcal{O}(1/\sqrt{n})$.

Concerning the second term, let us define
$$F(k) = \int_{k-1}^k \log \left|k-1/2-t \right|d\mu(t), \quad k = 1, \ldots, n.$$
Using the symmetry, notice that $F(n-k+1) = F(k)$. Therefore, we can suppose that $1 \leq k \leq (n+1)/2$. On the one hand, if $k = 1$, then
\begin{align*}
\abs{\frac{1}{\pi}\int_{0}^1 \frac{\log \abs{1/2 -t}}{\sqrt{t(n-t)}} dt} & \leq
\frac{1}{\pi}\int_{0}^1 \frac{\abs{\log \abs{1/2 -t}}}{\sqrt{t(n-t)}} dt \\
& \leq \frac{1}{\pi\sqrt{n-1}} \int_{0}^1 \frac{\abs{\log \abs{1/2 -t}}}{\sqrt{t}}dt =
\mathcal{O}\left(\dfrac{1}{\sqrt{n}}\right).
\end{align*}

On the other hand, if $2 \leq k \leq (n+1)/2$, then
\begin{align*}
\abs{\frac{1}{\pi} \int_{k-1}^k\frac{\log \abs{k-1/2-t}}{\sqrt{t(n-t)}}dt} & \leq
\frac{1}{\pi} \int_{k-1}^k\frac{\abs{\log \abs{k-1/2-t}}}{\sqrt{t(n-t)}}dt \\
& \leq \frac{1}{\pi\sqrt{(k-1)(n-k)}}\int_{k-1}^k\abs{\log \abs{k-1/2-t}}dt \\
& \leq \frac{1}{\pi\sqrt{(n-1)/2}}\int_0^1\abs{\log \abs{t-1/2}}dt
= \mathcal{O}\left(\dfrac{1}{\sqrt{n}}\right).
\end{align*}

Let us now estimate the third term. Fix $j = 1, \ldots, n$ with $j \neq k$, and notice that all $t \in [j-1,j]$ can be written as $t = j-1+t'$ for some $t' \in [0,1]$. Therefore,
$$\log\abs{k-j} - \log \abs{k-1/2-j+1-t'} = \log\abs{1+\dfrac{1-2t'}{2(k-j)}} = \mathcal{O}\left(\dfrac{1}{\abs{k-j}}\right),$$
where we have used that $\left|\frac{1-2t'}{2(k-j)}\right|\leq \frac12$.
Thus, there exists $c > 0$ such that
\begin{align*}
\sum_{\substack{j=1 \\ j \neq k}}^n  \int_{j-1}^j  \left(\log\abs{k-j} - \log \left|k-1/2 -t \right|\right)d\mu(t) & \leq \sum_{\substack{j=1 \\ j \neq k}}^n  \int_{j-1}^j  \dfrac{c}{\abs{k-j}}d\mu(t) \\
& = \sum_{\substack{j=1 \\ j \neq k}}^n  c\dfrac{\alpha_j}{\abs{k-j}} = \mathcal{O}\left(\dfrac{1}{\sqrt{n}}\right),
\end{align*}
where Lemma \ref{lem:properties_aj}.(b) has been used.

Combining these arguments, we conclude that $\abs{S} = \mathcal{O}(1/\sqrt{n})$. As a result, from the definition of $S$ and \eqref{eq:psigma}, we obtain that
\begin{equation}
\label{eq:potentialK}
\alpha_k I(\nu) + \sum_{\substack{j=1 \\ j\neq k}}^n \alpha_j \int_E \log|x-j-y| d\nu(y)= p_{\mu}(k-1/2) + S + \mathcal{O}\left(\dfrac{1}{\sqrt{n}}\right) = \log\left(\dfrac{n}{4}\right) + \mathcal{O}\left(\dfrac{1}{\sqrt{n}}\right).
\end{equation}
where \eqref{eq:properties_mu} has been used, what proves the claim and thus (a) and (b).

It is clear that (c) follows from (b).

(d) follows from similar ideas to the ones used to obtain (b). To see this, notice that $\abs{x-y} \leq 1/2$ for all $x,y \in E$, from which it follows that, for any given $k \in \N$,
\begin{equation}
\label{eq:log}
\log\abs{x-y-k}-\log(k) = \log\abs{1-\dfrac{x-y}{k}} = \mathcal{O}\left(\dfrac{1}{k}\right).
\end{equation}
Therefore, for every $x \in E$,
\begin{equation}
\label{eq:psigmaE}
p_{\sigma}(x) = \sum_{k = 1}^n\alpha_k\int_E\log\abs{x-y-k}d\nu(y) = \sum_{k = 1}^n\alpha_k\log(k) + \mathcal{O}\left(\dfrac{1}{\sqrt{n}}\right),
\end{equation}
where \eqref{eq:log} and Lemma \ref{lem:properties_aj}.(b) have been used.

From \eqref{eq:properties_mu}, $p_{\mu}(0)=\log(n/4)$ and we have 
\begin{align*}
S & := \sum_{k = 1}^n\alpha_k\log(k) - \log(n/4) = \sum_{k = 1}^n\alpha_k\log(k) - p_{\mu}(0) \\
& = \sum_{k = 1}^n\int_{k-1}^k\log(k)d\mu(t) - \sum_{k = 1}^n\int_{k-1}^k\log(t)d\mu(t)  = \sum_{k = 1}^n\int_{k-1}^k\log(k/t)d\mu(t).
\end{align*}

We claim that $\abs{S} = \mathcal{O}(1/\sqrt{n})$. To see this, notice that if $k = 1$, then
$$\int_{0}^1\log(1/t)d\mu(t) = \dfrac{1}{\pi}\int_0^1\dfrac{-\log(t)}{\sqrt{t(n-t)}}dt \leq \dfrac{1}{\pi\sqrt{n-1}}\int_0^1\dfrac{-\log(t)}{\sqrt{t}} dt = \mathcal{O}\left(\dfrac{1}{\sqrt{n}}\right).$$
On the other hand, if $k \geq 2$, notice that
$$\log\left(\dfrac{k}{k-1}\right) = \log\left(1+\dfrac{1}{k-1}\right) = \mathcal{O}\left(\dfrac{1}{k}\right).$$
Therefore, there exists some $c > 0$ such that the following holds:
\begin{align*}
\sum_{k=2}^n\int_{k-1}^k\log(k/t)d\mu(t) & \leq \sum_{k=2}^n\int_{k-1}^k\log(k/(k-1))d\mu(t) \\
& \leq \sum_{k=2}^n\int_{k-1}^k\dfrac{c}{k}d\mu(t) = c\sum_{k=2}^n\dfrac{\alpha_k}{k} = \mathcal{O}\left(\dfrac{1}{\sqrt{n}}\right),
\end{align*}
by Lemma \ref{lem:properties_aj}.(b). Notice that the claim follows from these arguments. As a result, from the definition of $S$ and \eqref{eq:psigmaE}, we obtain that
$$p_{\sigma}(x) = \log\left(\dfrac{n}{4}\right) + \mathcal{O}\left(\dfrac{1}{\sqrt{n}}\right).$$

Lastly, (e) follows from (b) and (d).
\end{proof}
\end{lemma}

At this point, we can proceed with the proof of Theorem \ref{thm:main-potential}.

\begin{proof}[\bf Proof of Theorem \ref{thm:main-potential}]
Consider the equilibrium measure $\beta$ for $K_n$. The definition of $\sigma$ depends on a finite translation of the equilibrium measure of $E$. Therefore, $\sigma (A)=0$ for every Borel subset $A$ with $\Capac(A)=0$. From Theorem \ref{thm:Frostman}, $I(\beta)=p_\beta (x)$ for nearly every $x\in K_n$. Moreover, it follows
$$\log(\Capac(K_n)) = \int_{K_n}p_{\beta}(x)d\sigma(x) = \int_{K_n}p_{\sigma}(x)d\beta(x),$$
where Fubini's Theorem has been used. Then, by Lemma \ref{lemma:sigma}.(c),
$$\abs{\log(\Capac(K_n)) - p_{\sigma}(y)} = \abs{\int_{K_n}(p_{\sigma}(x)-p_{\sigma}(y))d\beta(x)} \leq \int_{K_n}\abs{p_{\sigma}(x)-p_{\sigma}(y)}d\beta(x) = \mathcal{O}\left(\dfrac{1}{\sqrt{n}}\right)$$
for nearly every $y \in K_n$. Thus, the result follows from Lemma \ref{lemma:sigma}.(b) and \eqref{eq:properties_mu}.
\end{proof}

\section{Proof of Theorem \ref{thm:main-hardy}} \label{sec:proof2}
\begin{proof}[\bf Proof of Theorem \ref{thm:main-hardy}]Using Theorem \ref{thm:Nevanlinna}, we will restrict to the case where $\Omega$ is a regular Koenigs domain.

Consider $X = \C \setminus \Omega$, and define $B_N = \{z \in X : \abs{z} \leq N\}$ for all $N \in \N$. Notice that $\cup_{N \in \N}B_N = X$. Then, by Lemma \ref{lem:countable}, since $\Capac(X) > 0$, there must exists some $N \in \N$ such that $\Capac(B_N) > 0$. Moreover, given $z \in B_N$, define $E_z = \{w \in B_N : \abs{w-z} \leq 1/4\}$. Notice that $\cup_{z \in B_N}E_z = B_N$, which is a compact set. Then, there must exist a (finite) sequence $ \{x_n \}_{n=1,...,m}$ in $ B_N$, such that $\cup_{n =1}^{m}E_{x_n} = B_N$. In that case, using Lemma \ref{lem:countable} once more and the fact that $\Capac(B_N) > 0$, there must exist $x \in B_N$ such that $E_x \subseteq X$ with $\Capac(E_x) > 0$.

Without loss of generality, we can suppose that $x = 0$. By construction, $E_x \subseteq \overline{D(0,1/4)}$. With this in mind, define $E = -E_x$, and consider
$$\Omega^* = \C \setminus \bigcup_{n \in \N} (E+n).$$
Note that $\Omega \subseteq -\Omega^*$, which yields $\h(\Omega) \geq \h(-\Omega^*) = \h(\Omega^*)$. Therefore, it is enough to prove that $\h(\Omega^*) \geq 1/2$. To do so, Theorem \ref{thm:main-potential} may be used to estimate $\h(\Omega^*)$ in the following way:

As in Theorem \ref{thm:KimSugawa}, consider the harmonic measure $\omega_R(z) = \omega(z,F^*_R,\Omega^*_R)$, where $\Omega^*_R$ is the connected component of $\Omega^* \cap D(0,R)$ containing the origin and $F^*_R = \partial\Omega^*_R \cap \{w \in \C : |w| = R\}$. For a given $R > 0$, define $m = \max\{j \in \N : j \leq R/3\}$, and notice that
$$\Omega^*_R \cup \{w \in \C : m + 1/4 < \abs{w} < R\}) = D(0,R) \setminus K_{m}=: \Delta_R,$$
where $K_{m}:= \bigcup_{j = 1}^m(E+j)$. Therefore, consider $\beta_R(z) = \omega(z,F_R^*,\Delta_R)$ and $\tilde \beta_R(z) = \omega(z,\partial D(0,R),\Delta_R)$. By the maximum principle for harmonic functions \cite[Theorem 1.1.8]{Ransford} it is possible to show that $\omega_R(0) \leq \beta_R(0) \leq \tilde \beta_R(0)$. 

Furthermore, $\tilde \beta_R$ can be estimated through the use of potentials. To do so, we follow the ideas which were developed in Section \ref{sec:proof1}. Consider $\nu$ the equilibrium measure for $E$, and define the probability measure given by
$$\sigma(A) = \sum_{j = 1}^m\alpha_j\nu(A-j),
$$
for every Borel set $A \subset \C$. Notice that if $z \in \C$ is such that $\abs{z} = R$, then $p_{\sigma}(z) \geq \log(R-m-1/4) \geq \log(R/3)$ if $R$ is big enough. Take $\lambda_{m}=\inf\{p_{\sigma}(x):\, x\in K_{m}\}$. As $\sigma$ is supported by $K_{m}$, by the Minimum Principle for logarithmic potentials (see \cite[Theorem 3.1.4]{Ransford}), $\lambda_{m}=\inf\{p_{\sigma}(x):\, x\in \C\}$. Using both statements (a) and (b) of  Lemma \ref{lemma:sigma}, we have that 
$$\lambda_{m}=\log(m/4)+\mathcal{O}\left(\dfrac{1}{\sqrt{m}}\right) =p_{\sigma}(y)+ \mathcal{O}\left(\dfrac{1}{\sqrt{m}}\right) ,
$$ 
for nearly every $y \in K_m$. Therefore, the function $\gamma_{R}:=p_{\sigma}-\lambda_{m}$ is non-negative on $\C$, harmonic in $\C\setminus K_{m}$, and for $z \in \C$ such that $\abs{z} = R$ we have
\begin{equation*}
\begin{split}
\gamma_{R}(z)&\geq \log(R/3)-\lambda_{m}=\log(R/3)-\log(m/4)+\mathcal{O}\left(\dfrac{1}{\sqrt{m}}\right)\\
&\geq \log(R/3)-\log(R/12)+\mathcal{O}\left(\dfrac{1}{\sqrt{m}}\right)=\log(4)+\mathcal{O}\left(\dfrac{1}{\sqrt{m}}\right).
\end{split}
\end{equation*}
Therefore, $\gamma_{R}(z)\geq 1$ whenever $|z|=R$ for $m$ large enough.
Thus $\tilde \beta_R(0)\leq \gamma_{R}(0)=p_{\sigma}(0)-\lambda_{m}$, for $m$ (and then $R$) large enough.

By Lemma \ref{lemma:sigma}.(e),  for nearly every $y \in K_m$,
$$
\gamma_{R}(0)=p_{\sigma}(0)-\lambda_{m}=p_{\sigma}(0)-p_{\sigma}(y)+ \mathcal{O}\left(\dfrac{1}{\sqrt{m}}\right) =\mathcal{O}\left(\dfrac{1}{\sqrt{m}}\right) =\mathcal{O}\left(\dfrac{1}{\sqrt{R}}\right).
$$

With this in mind, Theorem \ref{thm:KimSugawa} yields
$$\h(\Omega^*) = \liminf_{R \to + \infty}\left(-\dfrac{\log\omega_R(0)}{\log R}\right) \geq \liminf_{R \to + \infty}\left(-\dfrac{\log\gamma_R(0)}{\log R}\right) \geq \dfrac{1}{2}.$$
\end{proof}

\end{document}